\documentclass{amsart}
\usepackage{amssymb}
\usepackage{amscd}
\usepackage{amsmath}
\usepackage{enumerate}
\usepackage[dvips]{graphicx}
\newtheorem{theo}{Theorem}
\newtheorem{lema}[theo]{Lemma}

\newtheorem{definition}[theo]{Definition}

\newtheorem{remark}[theo]{Remark}

\newtheorem{problem}[theo]{Problem}

\usepackage{tikz}

\newcommand{\CC}{{\mathbb{C}}}

\newcommand{\QQ}{{\mathbb{Q}}}

\newcommand{\SSS}{{\mathbb{S}}}
\newcommand{\ZZ}{{\mathbb{Z}}}

\newcommand{\calO}{{\mathcal{O}}}

\newcommand{\calU}{{\mathcal{U}}}

\newcommand{\comp}{{\circ}}

\newcommand{\Sing}{{\mathrm{Sing}}}
\begin{document}
\title[The Nash problem]{The Nash problem from a geometric and topological perspective}
\author{J. Fern\'andez de Bobadilla}
\address{BCAM-Ikerbasque. Alameda de Mazarredo 14. 48009 Bilbao}
\email{jbobadilla@bcamath.org}
\author{M. Pe Pereira}
\address{Facultad de Ciencias Matematicas- UCM, Plaza de Ciencias, 3,  Ciudad Universitaria, 28040 MADRID} 
\email{maria.pe@mat.ucm.es}
\thanks{First author is supported by ERCEA 615655 NMST Consolidator Grant, MINECO by the project reference 
 MTM2016-76868-C2-1-P, by the Basque Government through the BERC 2014-2017 program, 
by Spanish Ministry of Economy and Competitiveness MINECO:  MTM2016-76868-C2-1-P, BCAM Severo Ochoa excellence accreditation SEV-2013-0323 and 
by Bolsa Pesquisador Visitante Especial (PVE) - Ciências sem Fronteiras/CNPq Project number:  401947/2013-0. Part of this work was carried 
when he was a member at IAS (Princeton); he thanks AMIAS for support. Second author is supported by  MINECO by the project reference 
 MTM2016-76868-C2-1-P, by Bolsa Pesquisador Visitante Especial (PVE) - Ciências sem Fronteiras/CNPq Project number:  401947/2013-0 and by 
NMST Consolidator Grant.}

\date{14-1-2018}
\begin{abstract}
We survey the proof of the Nash conjecture for surfaces and show how geometric and topological ideas developed in previous articles by the authors influenced it. Later we summarize the main ideas in the higher dimensional statement and proof by de Fernex and Docampo. We end the 
paper by explaining later developments on generalized Nash problem and on Koll\'ar and Nemethi holomorphic arcs. 
\end{abstract}


\maketitle

\section{Introduction}

The Nash problem~\cite{Na} was formulated in the sixties (but published later) in the attempt to understand the relation between the 
structure of resolution of singularities of an algebraic variety $X$ over a field of characteristic $0$ and the space of arcs 
(germs of parametrized curves) in the variety.
He proved that the space of arcs centered at the singular locus (endowed with an infinite-dimensional algebraic variety structure) 
has finitely many irreducible components and proposed to study the relation of these components with the 
essential irreducible components of the exceptional set of a resolution of singularities. 

An irreducible component $E_i$ of the exceptional divisor of a
resolution of singularities is called \emph{essential}, if given any other resolution 
the birational transform of $E_i$ to the second resolution is an irreducible component of the exceptional divisor. Nash defined a mapping from the set of
irreducible components of the
space of arcs centered at the singular locus to the set of essential components of a resolution as follows: he assigns to each component $W$ 
of the space of arcs centered at the singular locus the unique component of the exceptional set
which meets the lifting of a generic arc of $W$ to the resolution. Nash established the injectivity of this mapping. 
For the case of surfaces it seemed plausible for him that the mapping is also surjective, and posed
the problem as an open question. He also proposed to study the mapping in the higher dimensional case.
Nash resolved the question positively for the surface $A_k$ singularities, and in analyzing the higher dimensional $A_k$ singularities he could not 
prove the bijectivity for $A_4$. 

As a general reference for the Nash problem the reader may look at~\cite{Na} and~\cite{IK}.

Bijectivity of the Nash mapping was shown for many classes of surfaces 
(see \cite{Go},\cite{IK},\cite{I1},\cite{I2},\cite{Le},\cite{LR1},\cite{Mo},\cite{Pe},\cite{Pl},\cite{Pl2},\cite{PP1},\cite{PlSp},\cite{Re1},\cite{Re2}).
The techniques leading to the proof of each of these cases are 
different in nature, and the proofs are often complicated. It is worth to notice that even for the case of
the rational double points not solved by Nash a complete proof has to be awaited
until 2011: see~\cite{Pe}, where the problem is solved for any quotient surface singularity, and also \cite{Pl2} and~\cite{PlSp} for the cases of $D_n$ and $E_6$.
In~\cite{Bo} it is shown 
that the Nash problem for surfaces only depends on the topological type of the singularity. In 2012 the authors of the present paper established in the 
afirmative the Nash question for the general surface case~\cite{annals}. The proof we found was of a topological nature, and it is essential to work with
convergent arcs and their convergent deformations. This motivated J. Koll\'ar and A. Nemethi to pursue the study of convergent arcs and deformations 
in~\cite{KoNe}. The topological ideas of~\cite{Bo} and \cite{Pe} also had an impact on the generalized Nash problem; in~\cite{BoPP} P. Popescu-Pampu and the authors
of the present survey show that the generalized Nash problem is of topological nature and explore the relation and applications of this problem to
Arnol'd classical adyacency problem.

It is well known that birational geometry of surfaces is much simpler than in higher dimension. This fact reflects on the Nash problem: 
Ishii and Koll\'ar showed in~\cite{IK} a 4-dimensional example with a non-bijective Nash mapping. In the same paper they showed the bijectivity of the Nash 
mapping for toric singularities of arbitrary dimension. Other advances in the higher dimensional case include~\cite{PP2},~\cite{Go},~\cite{LR}.
In 2013 T. de Fernex~\cite{deF} found the first counter-examples to the Nash question; further counterexamples, and a deeper understanding of how they appear
was provided by J. Johnson and J. Koll\'ar in~\cite{JK}. There it was proved that the threefold $A_4$ 
$$x^2+y^2+z^2+w^5=0,$$
the example that Nash left unfinished, was indeed a counter-example! In 2016 T. de Fernex and R. Docampo~\cite{deFDo} 
proved that terminal divisors are at the image of the Nash map. Since at the surface case, terminal and essential divisors are precisely the same, this 
seems to be the correct higher dimensional generalization. It would remain, however, to characterize which essential non-terminal divisors at at the image
of the Nash map. 

For other modern review articles concerning the Nash problem the reader may consult~\cite{PLSp2},~\cite{Ko},~\cite{JK2}  

In this paper we explain how geometric and topological techniques contributed to the development of the proof of the Nash conjecture, and how they relate
with other viewpoints and further developments. 

In Sections 2-5 explain our proof of the $2$ dimensonal case in a non-technical way, pointing to the main new ideas appearing in it.
We present a proof for the case in which the minimal resolution has a strict normal crossings exceptional divisor.
In this case all new essential ideas 
already appear, but the ammount of technicalities can be reduced drastically.
We include enough pictures so thar the reader can grasp what is going on
in an easy and intituitive way. 

In Section 6 we emphasize the notion of returns, which was discovered in~\cite{Pe} and was crucial for the development of the general proof. 
We also take the opportunity to comment on deformation techniques that were useful to establish the hard cases of $E_6$, $E_7$ and $E_8$.

In Section 7  we explain the relation of our proof with the higher dimensional one of \cite{deFDo}. We do it by giving a short exposition of their 
proof that we believe condense all main ideas.

In Section 8 we summarize our contribution with P. Popescu-Pampu on the generalized Nash problem~\cite{BoPP}, 
and its impact on Arnol'd classical adyacency problem. Here we use the techniques of~\cite{Bo} to show that the generalized Nash problem is of 
topological nature.

Finally in Section 9 we explain the relation of our ideas with the further developments
of more geometric-topological nature by Koll\'ar and Nemethi~\cite{KoNe}.

\section{The idea of the proof for surfaces}

Let $(X,O)$ be a surface singularity defined over an algebraically closed field of $0$ characteristic. Let 
$$\pi:(\tilde{X},E)\to (X,O)$$
be the minimal resolution of singularities, which is an isomorphism outside the exceptional divisor $E:=\pi^{-1}(O)$.
Consider the decomposition 
$E=\cup_{i=0}^r E_i$ of $E$ into irreducible components. These irreducible components are the essential components of $(X,O)$.

Given any irreducible component $E_i$ we denote by $N_{E_i}$ the Zariski closure in the arc space of $X$ of the set of 
non-constant arcs whose lifting to the resolution is centered at $E_i$. These Zariski closed subsets are irreducible and each irreducible component
of the space of arcs is equal to some $N_{E_i}$ for a certain component $E_i$. 
The Nash mapping is the map assigning to each irreducible component $N_{E_i}$ the exceptional divisor $E_i$. Injectivity is immediate. The 
Nash problem is about determining whether the Nash mapping is bijective. 

The Nash mapping is not bijective if and only if there exist two different irreducible components
$E_i$ and $E_j$ of the exceptional divisor of the minimal resolution such that we have the inclussion 
$N_{E_i}\subset N_{E_j}$ (see \cite{Na}). Such inclusions were called {\em adjacencies} in~\cite{Bo}.

An application of the Lefschetz principle allows to reduce to the case in which the base field is $\CC$. Details are at~\cite{annals}.
We make this assumption for the rest of the paper.
Moreover, the case of a non-normal surface follows from the normal surface case easily (see Section 6 in \cite{annals}). 
Then, we assume  $(X,O)$ to be a complex normal surface singularity.

The idea of the proof is as follows. We reason by contradiction.
Let $(X,O)$ be a normal surface singularity and
$$\pi:\tilde{X}\to (X,O)$$
be the minimal resolution of singularities. Assume that the Nash mapping is not bijective.
Then, by a Theorem of~\cite{Bo} there exists a convergent wedge
$$\alpha:(\CC^2,O)\to (X,O)$$
with certain precise properties (see Definition~\ref{wedadj}). As in~\cite{Pe}, taking a suitable representative we may view $\alpha$ as a
uniparametric family of mappings
$$\alpha_s:\calU_s\to (X,O)$$
 from a family of domains $\calU_s$ to $X$ with the property that each $\calU_s$ is diffeomorphic to a disk.
For any $s$ we consider the lifting
$$\tilde{\alpha}_s:\calU_s\to\tilde{X}$$
to the resolution. Notice that $\tilde{\alpha}_s$ is the normalization mapping of the image curve. 

On the other hand, if we denote by $Y_s$ the image of $\tilde{\alpha}_s$ for $s\neq 0$, then we may consider the
limit divisor $Y_0$ in $\tilde{X}$ when $s$ approaches $0$. This limit divisor consists of the union of the image of 
$\tilde{\alpha}_0$ and certain components of the exceptional divisor of the resolution whose multiplicities are easy to compute. We prove an upper
bound for the Euler characteristic of the normalization of any reduced deformation of $Y_0$ in terms of the following data:
the topology of $Y_0$, the multiplicities of its components and
the set of intersection points of $Y_0$ with the generic member $Y_s$ of the deformation. Using this bound we show that the Euler characteristic of the 
normalization of $Y_s$ is strictly smaller than one. This contradicts the fact that the normalization is a disk.

In the next three sections we fill the details of the sketch above.

\section{Turning the problem into a problem of convergent wedges}

The germ $(X,O)$ is embedded in an ambient space $\CC^N$. Denote by $B_{\epsilon}$ the closed ball of radius $\epsilon$ centered at
the origin and by $\SSS_{\epsilon}$ its boundary sphere. Take a \emph{Milnor radius} $\epsilon_0$ for $(X,O)$ in $\CC^N$,
that is, we choose $\epsilon_0>0$ such that for a certain representative $X$ and any radius $0<\epsilon\leq \epsilon_0$ we have that
all the spheres $\mathbb{S}_{\epsilon}$  are transverse to $X$ and $X\cap \SSS_{\epsilon}$ is a closed subset of $\SSS_{\epsilon}$
(see~\cite{Mil} for a proof of its existence). In particular $X\cap B_{\epsilon_0}$ has conical structure. 
From now on we will denote by $X_{\epsilon_0}$ the \emph{Milnor representative} $X\cap B_{\epsilon_0}$ and by 
$\tilde{X}_{\epsilon_0}$ the resolution of singularities $\pi^{-1}(X_{\epsilon_0})$.

We recall some terminology and results from~\cite{Bo}. Consider coordinates $(t,s)$ in the germ $(\CC^2,O)$. A \emph{convergent wedge} is a 
complex analytic germ
$$\alpha:(\CC^2,O)\to (X,O)$$
which sends the line $V(t)$ to the origin $O$. Given a wedge $\alpha$ and a parameter value $s$, the arc
$$\alpha_s:(\CC,0)\to (X,O)$$
is defined by $\alpha_s(t)=\alpha(t,s)$. The arc $\alpha_0$ is called {\em the special arc} of the wedge. For small enough $s\neq 0$
the arcs $\alpha_s$ are called \emph{generic arcs}.

Any non-constant arc
$$\gamma:(\CC,0)\to (X,O)$$
admits a unique lifting  to $(\tilde{X},O)$ that we denote by $\tilde{\gamma}$.

\begin{definition}[\cite{Bo}]
\label{wedadj}
A convergent wedge $\alpha$ {\em realizes an adjacency} $N_{E_i}\subset N_{E_j}$ (with $j\neq i$) if and only if the lifting $\tilde{\alpha}_0$ of the special 
arc meets $E_i$ transversely at a non-singular point of $E$ and the lifting $\tilde{\alpha}_s$ of a generic arc satisfies $\tilde{\alpha}_s(0)\in E_j$.
\end{definition}

Our proof is based on the following Theorem, which is the implication ''$(1)~\Rightarrow~(a)$~'' of Corollary~B~of~\cite{Bo}:
\begin{theo}[\cite{Bo}]
\label{nashwedges}
An essential divisor $E_i$ is in the image of the Nash mapping if there is no other essential divisor 
$E_j\neq E_i$ such that there exists a convergent wedge realizing the adjacency $N_{E_i}\subset N_{E_j}$.
\end{theo}

The proof in \cite{Bo} of this theorem has two parts. The first consists of proving that if there is an adjacency then there exists 
a {\em formal wedge}
$$\alpha:Spec(\CC[[t,s]])\to (X,O)$$
realising the adjacency. For that, firstly it is used a Theorem of A. Reguera~\cite{Re} which produces 
wedges defined over large fields. Then a specialisation argument is performed to produce a wedge defined over the
base field $\CC$. This was done independently in~\cite{LR}. The second part is an argument based on D. Popescu's 
Approximation Theorem, which produces the convergent wedge from the formal one. 

Then, to prove the Nash Conjecture we reason by contradiction and  by Theorem \ref{nashwedges} we assume that there exists a convergent wedge $\alpha:(\CC^2,O)\to (X,O)$ that realises some adjacency  $N_{E_0}\subset N_{E_j}$. 

\section{Reduction to an Euler characteristic estimate}
\label{representatives}

Following~\cite{Pe} we shall work with representatives rather than germs in order to get richer information about the geometry of the possible wedges. 

Shrinking $\epsilon$ is neccessary, we can choose a Milnor representative of  $\alpha_0$, say $\alpha_0|_{U}$, with $U$ diffeomorphic to a disk, 
such that $\alpha_0|_U^{-1}(\partial X_\epsilon)=\partial U$ and the mapping $\alpha_0|_U$ is transverse to any sphere $\SSS_{\epsilon'}$ for any
$0<\epsilon'\leq \epsilon$.

Moreover, we can consider $U$ 
such that for a positive and small enough $\delta$ the mapping $\alpha$ is defined in 
$U\times D_\delta$. Note also that we can assume $\alpha_0|_U$ injective and consecuently $\alpha_s|_U$ generically one-to-one for $s$ small enough (see \cite{annals} for details).

We consider the mapping
$$\beta:(\CC^2,(0,0))\to(\CC^N\times\CC,(O,0))$$
given by $\beta(t,s):=(\alpha(t,s),s)$ 
and its restriction
$$\beta|_{U\times D_\delta}:U\times D_\delta\to X\times D_\delta.$$
We denote by $pr$ the projection of $U\times D_\delta$ onto the second factor. 

\begin{figure}
\includegraphics[width=60mm]{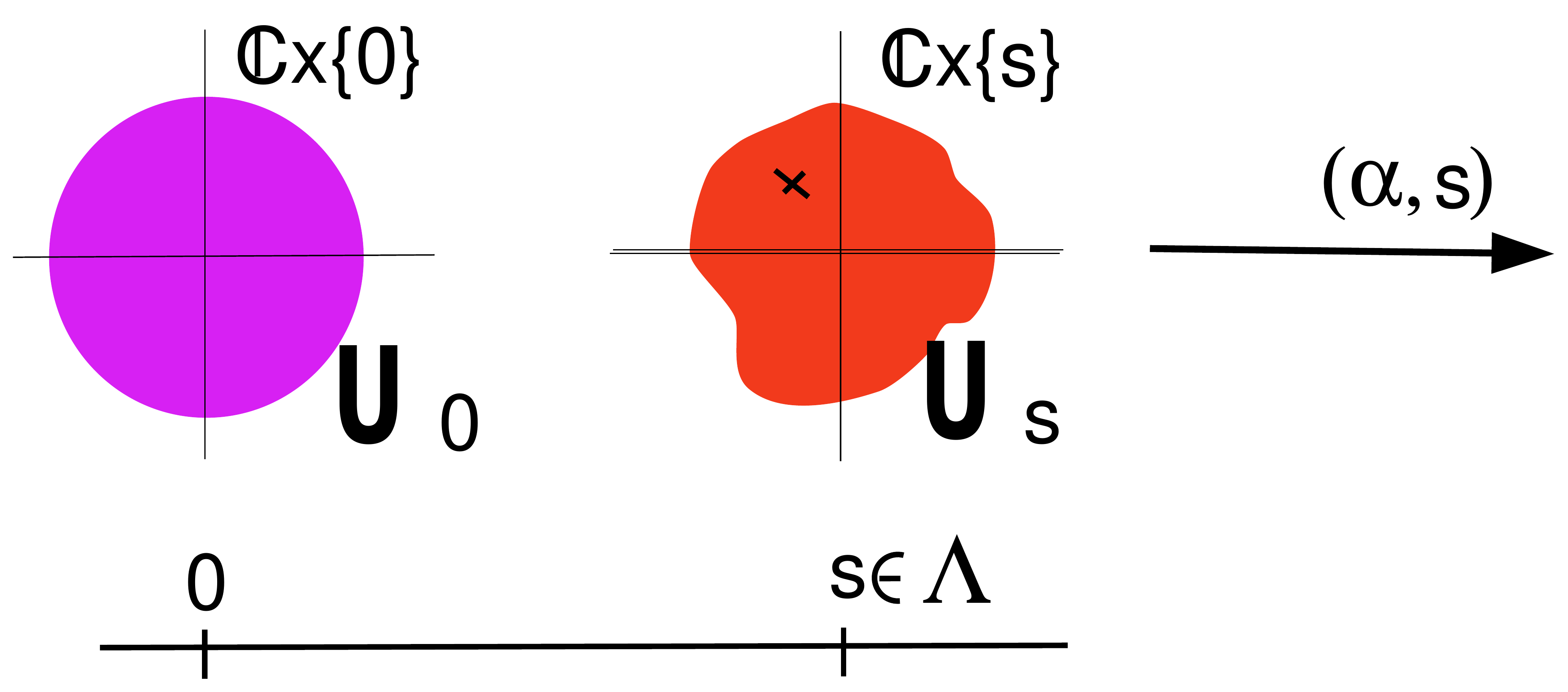}
\includegraphics[width=60mm]{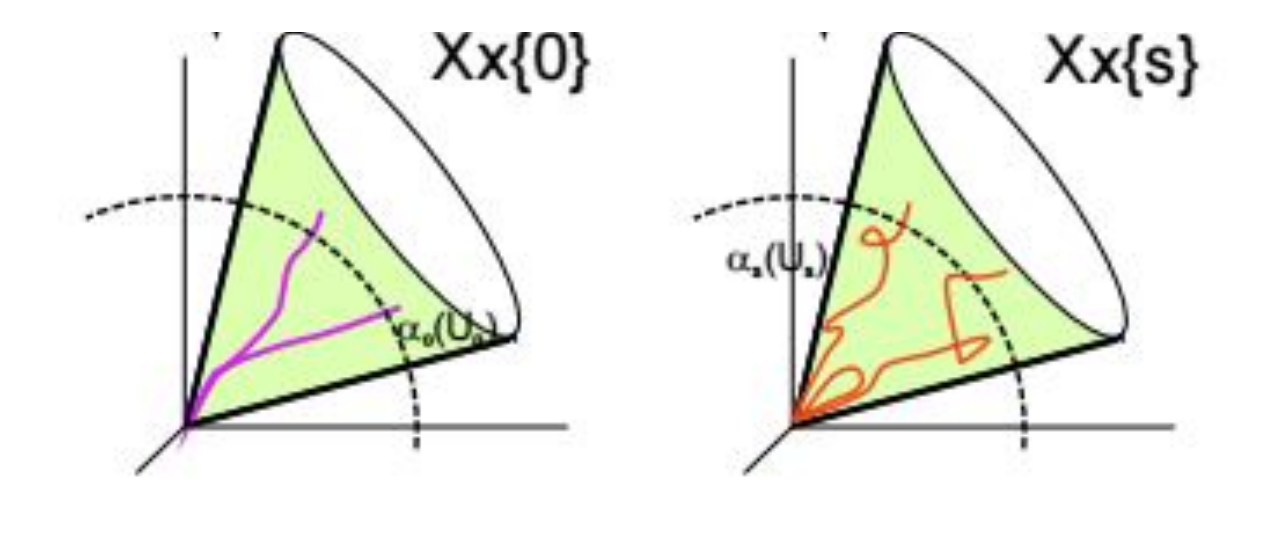}
\caption{A wedge representative $\alpha:\calU\to X\times \Lambda$ and the  representives $\alpha_0|_{\calU_0}$ and $\alpha_s|_{\calU_s}$.}
\end{figure}

The following Lemma is proved using transversality arguments, together with Ehresmann Fibration Theorem, and the method is nowadays classical in
Singularity Theory. Details are at~\cite{annals}.

\begin{lema}\label{lem:wedge_repre}
After possibly shrinking $\delta$, we have that there exists $\epsilon>0$ such that, defining 
$$\calU:=\beta|_{U\times {D}_\delta}^{-1}(X_\epsilon\times{D}_\delta)$$
we have that
\begin{enumerate}[(a)]
\item the restriction $\beta|_{\calU}:{\calU}\to {X}_{\epsilon}\times {D}_\delta$ is a proper
and finite morphism of analytic spaces, 
\item the set $\beta({\calU})$ is a two-dimensional closed 
analytic subset of ${X}_{\epsilon}\times{D}_\delta$,
\item for any $s\in D_\delta$ the restriction $\beta|_{U\times\{s\}}$ is transverse to $\SSS_\epsilon\times\dot{D}_\delta$,
\item the set $\calU$ is a smooth manifold with boundary $\beta|_{\calU}^{-1}(\partial X_\epsilon\times {D}_\delta)$ and
\item for any $s\in D_\delta$ the intersection $\calU\cap (\CC\times\{s\})$ is diffeomorphic to a disk.
\end{enumerate}
\end{lema}

We will denote by $\calU_s$ the fibre $pr|_{\calU}^{-1}(s)$. The fact that every $\calU_s$ is a disk is a key in the proof as it was in the final step of the proof of the main result of~\cite{Pe}.


 
Now we consider the image $H:=\beta({\calU})$. For every $s\in D_\delta$ the fibre $H_s$, by the natural projection onto $D_\delta$, is the image of the representative
$$\alpha_s|_{{\calU}_s}:{\calU}_s\to {X}_\epsilon.$$

\begin{figure}\label{fig:lift_wedge}
\begin{center}
      \includegraphics[width=40mm]{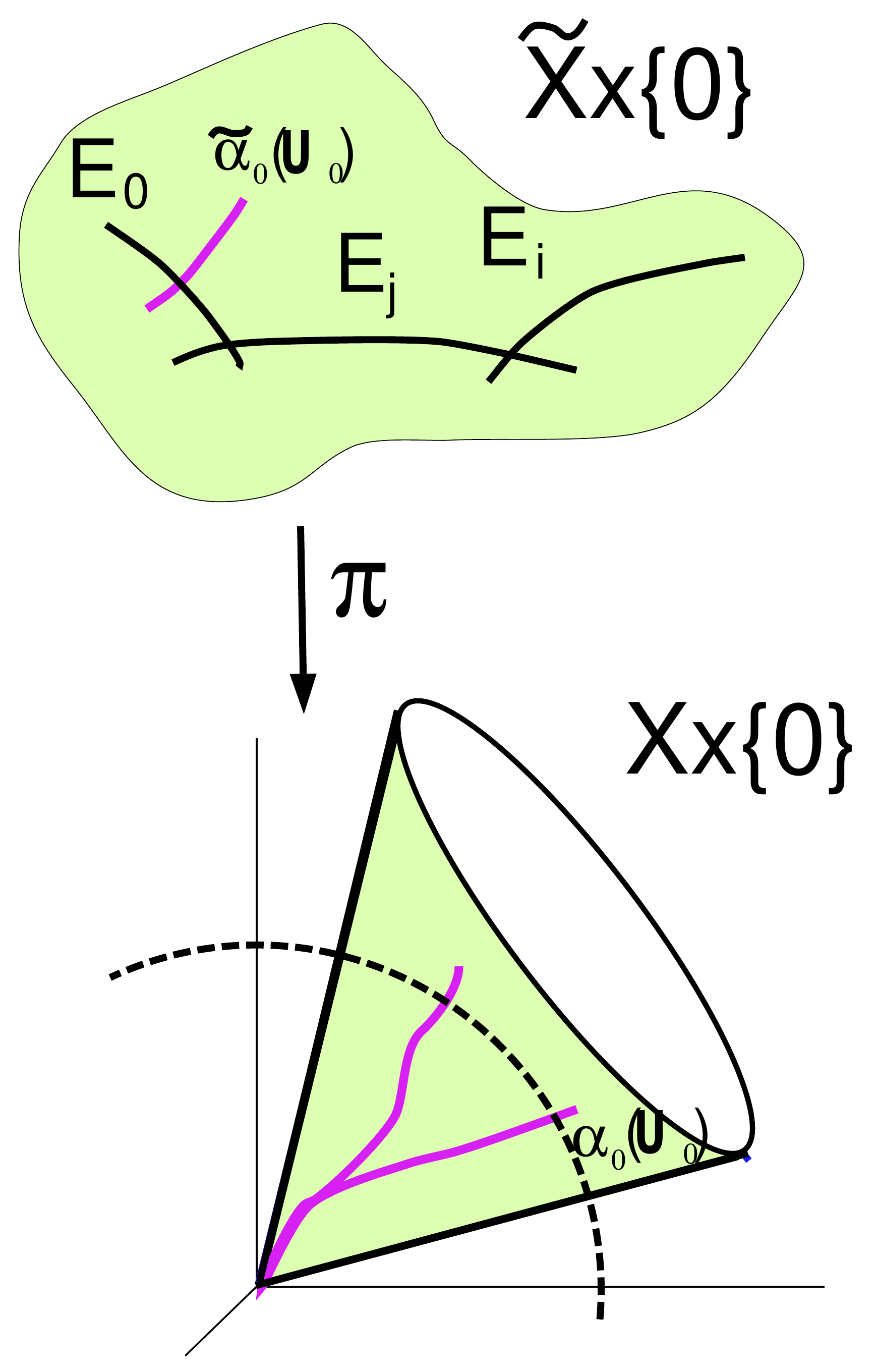}
      \includegraphics[width=40mm]{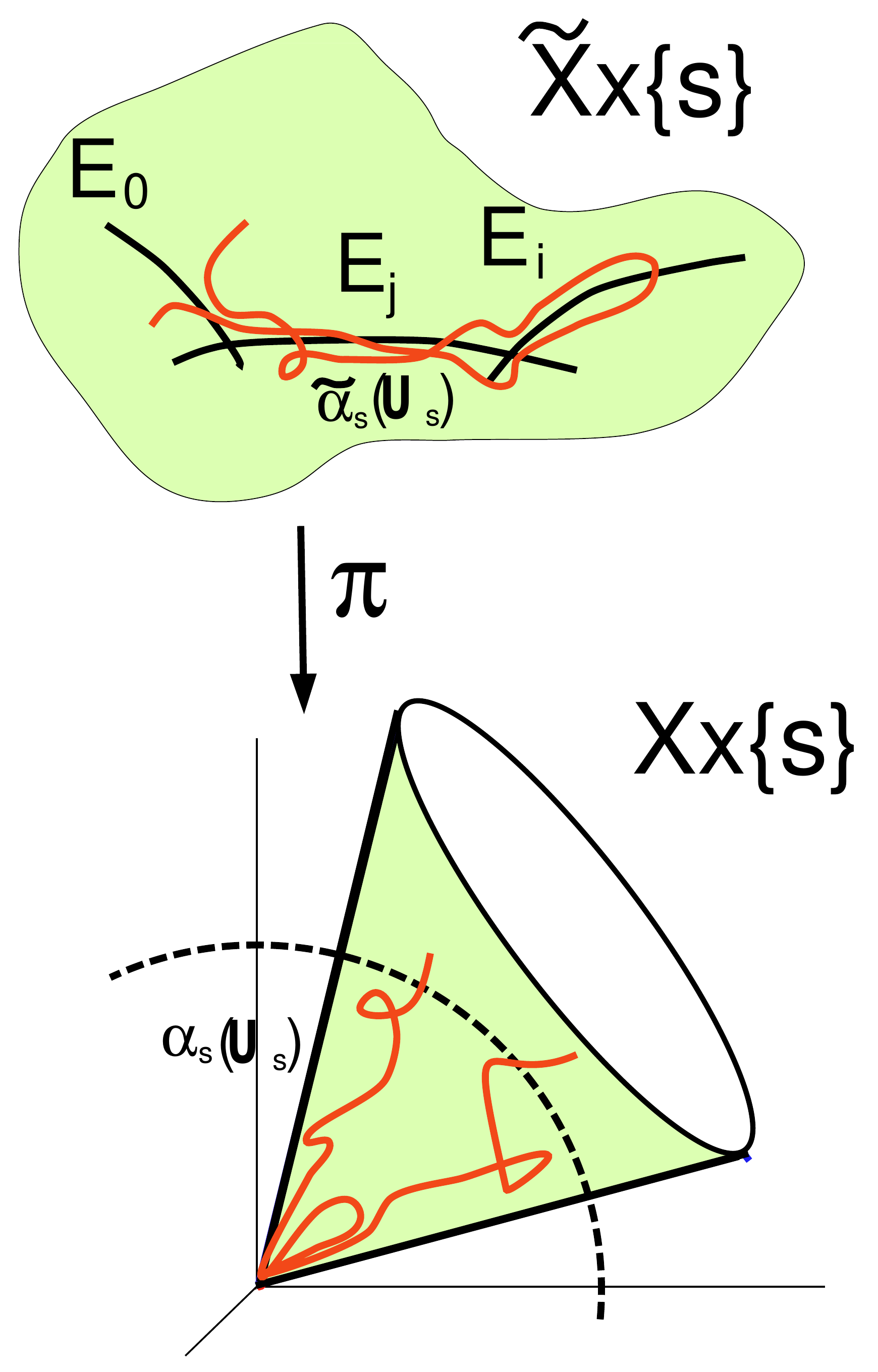}
      \caption{}
			\end{center}	
\end{figure}

Given the minimal resolution of singularities 
$$\pi:\tilde{X}_\epsilon\to X_\epsilon$$ we consider the mapping
$$\sigma:\tilde{X}_\epsilon\times D_\delta\to X_\epsilon\times D_\delta$$
defined by $\sigma(x,s)=(\pi(x),s)$. Note that the mapping $\sigma$ is an isomorphism outside $E\times D_\delta$. We denote by $Y$ the strict transform of $H$ by $\sigma$ in $\tilde{X}_\epsilon\times D_\delta$, that is the analytic Zariski closure in $\tilde{X}_\epsilon\times D_\delta$ of
\begin{equation}
\label{eq:Y}
\sigma^{-1}(H\setminus (\{O\}\times D_\delta)).
\end{equation}
The space (\ref{eq:Y}) is an irreducible surface, thus so is its closure $Y$. Since $\tilde{X}_\epsilon\times D_\delta$
is a smooth threefold, the surface $Y$ considered with its reduced structure is a Cartier divisor (that is, a 
codimension one analytic subset whose sheaf of ideals is locally principal). 
We denote by $Y_s$ the intersection $Y\cap(\tilde{X}\times\{s\})$.

The indeterminacy locus of the mapping $\sigma^{-1}\comp\beta|_{\calU}$ has codimension $2$, hence 
reducing $\epsilon$ and
$\delta$ if necessary, we can assume that the origin $(0,0)\in \calU$ is the only indeterminacy point. Denote by
$$\tilde{\beta}:\calU\setminus\{(0,0)\}\to\tilde{X}_\epsilon\times D_\delta$$
the restriction of $\sigma^{-1}\comp\beta|_{\calU}$ to its domain of definition $\calU\setminus\{(0,0)\}$.
Observe that we have the equality
$$\tilde{\beta}(\calU\setminus\beta^{-1}(\{O\}\times D_\delta))=\sigma^{-1}(H\setminus (\{O\}\times D_\delta)).$$
Consequently $Y$ is the analytic Zariski closure of 
$\tilde{\beta}(\calU\setminus\{(0,0)\})$
and moreover we have the equality
\begin{equation}
\label{pijotada}
Y\cap(\tilde{X}_\epsilon\times(D_\delta\setminus\{0\}))=\tilde{\beta}(\calU\setminus\calU_0).
\end{equation}

For any $s\in D_\delta$ there exists a unique lifting
$$\tilde{\alpha}_s:\calU_s\to\tilde{X}_\epsilon$$
such that $\alpha_s=\pi\comp\tilde{\alpha}_s$. Obviously for $s\neq 0$ we have the equality 
$\tilde{\beta}(t)=(\tilde{\alpha}_s(t),s)$ for any $t\in\calU_s$. This, together with Equality~(\ref{pijotada}), implies
the equality
\begin{equation}\label{eq:Ys}Y_s=\tilde{\alpha}_s(\calU_s).\end{equation}

Since $Y$ is reduced, perhaps shrinking $\delta$, we can assume that $Y_s$ is reduced. Since $\alpha_s$ is proper and 
generically one to one, and $\calU_s$ is smooth, we have that the mapping
$$\tilde{\alpha}_s:\calU_s\to Y_s$$
is the normalisation of $Y_s$. 
%

Now we describe the divisor $Y_0$. 
It is clear that all the components, except $\tilde{\alpha}_0(\calU_0)$ live above the origin that is the only indeterminacy point of $\tilde{\beta}$. That is, the divisor $Y_0$ decomposes as a sum
\begin{equation}
\label{decomposition}
Y_0=Z_0+\sum_{i=0}^ra_iE_i.
\end{equation}
where we have denoted 
$Z_0:=\tilde{\alpha}_0(\calU_0)$. 
This divisor $Y_0$ has the following properties: 
\begin{itemize}
\item[(i)] it is reduced at $Z_0\setminus E$ since $\sigma$ is an isomorphism outside
$E\times D_\delta$ and $H_0$ is reduced out of the origin. 
\item[(ii)] $Z_0$ intersects transversaly $E_0$ in a smooth point 
\item[(iii)]  all the $a_i$'s are non-negative since the divisor $Y_0$ is effective
\item[(iv)] some $a_i$ are strictly non zero, in particular $a_0$, since $\alpha$ realises an adjacency and then $\tilde{\beta}|_\calU$ has indeterminacy.
\end{itemize}

Assuming the existence of a wedge realizing an adjacency we have found a deformation $Y_s$ of some $Y_0$ (as in (\ref{decomposition}) and satifying (i)-(iv)) that has the following properties: 
\begin{itemize}
\item[(a)]  $Y_s$ is reduced for $s\neq 0$ small enough, 
\item[(b)] its normalisation, that is $\calU_s$, is diffeomorphic to a disk 
\item[(c)] its boundary, that is $\tilde{\alpha}_s(\partial \calU_s)$, is an $\SSS^1$ that degenerates to the boundary of $Y_0$, that is $\tilde{\alpha}_0(U)\cap \partial \tilde{X}_\epsilon$. 
\item[(d)] $Y_s$ meets $E_j\neq 0$.
\end{itemize}
The remaing part of the proof consists in proving that the euler characteristic of the normalization of such a deformation $Y_s$ of $Y_0$ is less or equal than $0$ which contradicts (b).

\section{The Euler characteristic estimates.}

To simplify the computation of the euler characteristic estimates we assume the minimal resolution of  $(X,O)$ has as exceptional divisor a simple normal 
crossings divisor. This is the first case that we discovered. The general case is technically more elaborate, but follows essentially the same ideas. 
It may be checked in~\cite{annals}. 

Let $Y_0$ be a Cartier divisor as in (\ref{decomposition}) that satisfies (i)-(iv). Consider a deformation $Y_s$ of $Y_0$ satisfying (a)-(d). 
Let $n:\calU_s\to Y_s$ be its normalization. 

\begin{figure}\label{fig:tubular}[h!]\begin{center}
      \includegraphics[width=100mm]{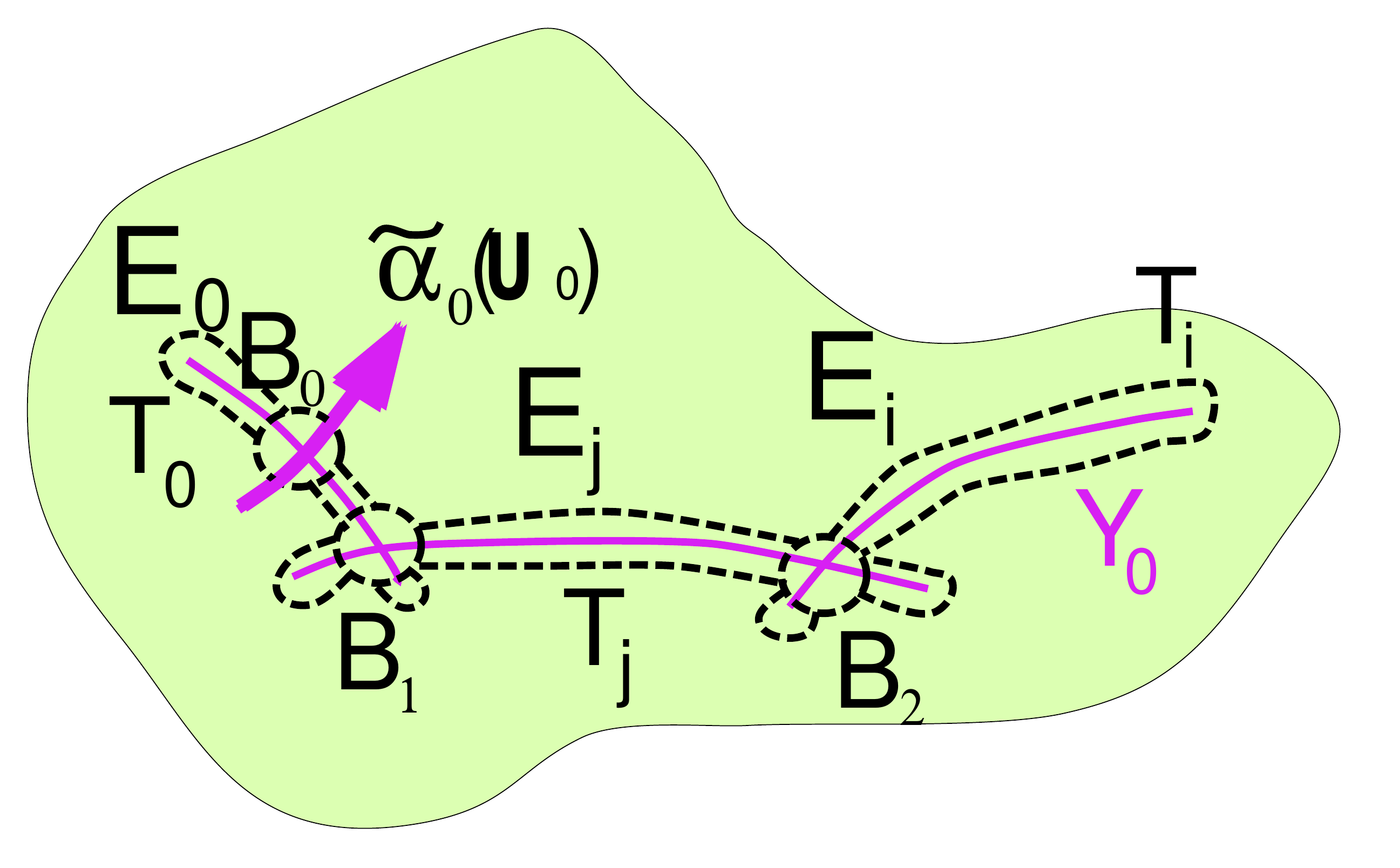}
      \caption{Adapted tubular neighbourhood of $Y_0$}
\end{center}	
\end{figure}

\begin{figure}\label{fig:disk}\begin{center}
  \includegraphics[width=40mm]{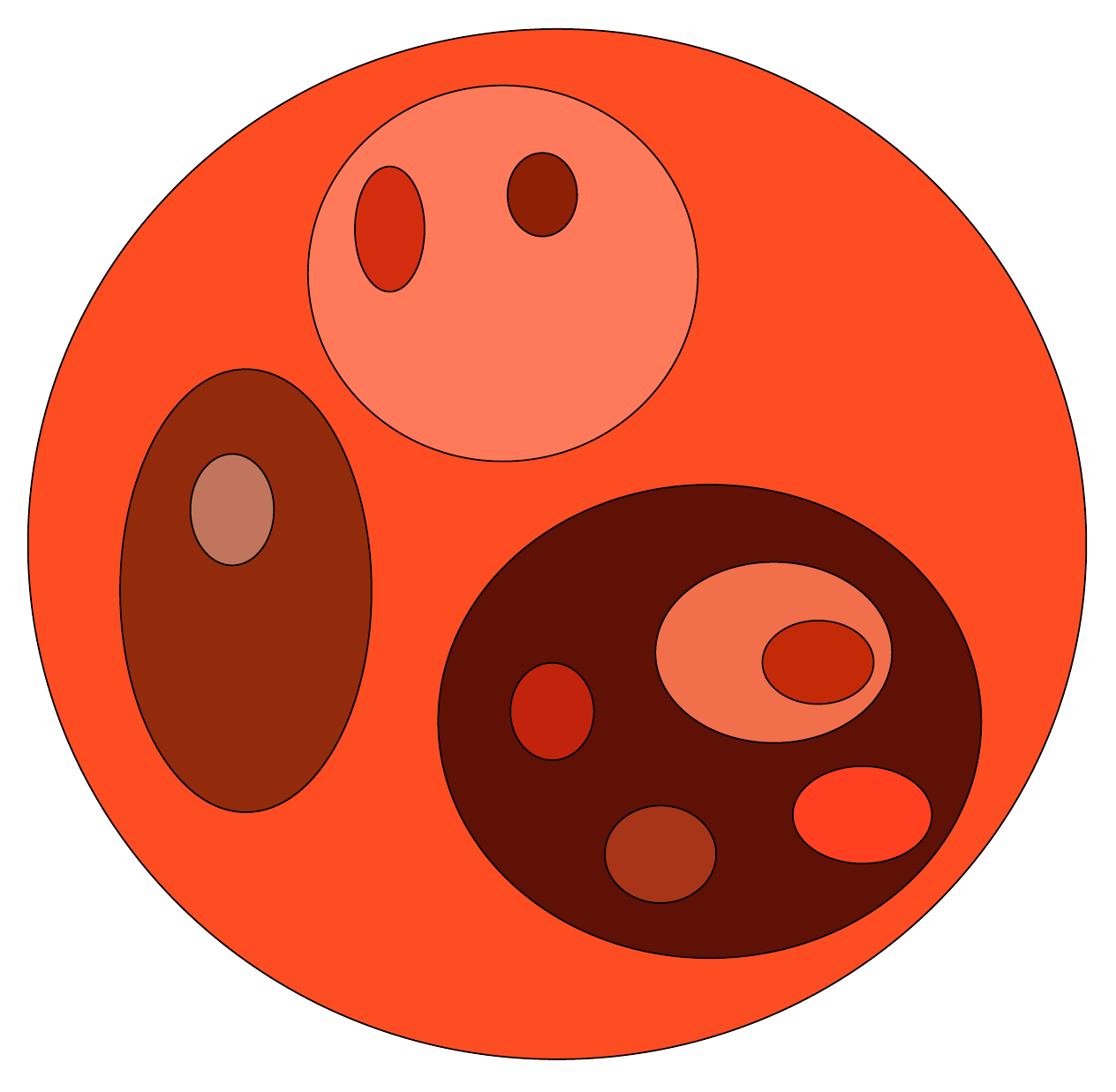}
			\end{center}	
\caption{The normalization $\calU_s$ of $Y_s$ inside the tubular neighbourhood of $Y_0$ is a disc. In the picture we see the result of cutting $Y_s$ along the boundary of the milnor balls $B_i$ around the normal crossings of $Y_0$. Each piece is either $n^{-1}(B_i)$ or $n^{-1}(T_j)$. The exterior piece that we call $A$  satisfies that $n(A)$ is contained in $B_0$.}
\end{figure}

We consider a tubular neighbourhood of $Y_0$ inside $\tilde{X}$ as a union of the following sets: 

\begin{itemize}
\item a Milnor ball $B_0:=B(Z_0\cap E, \epsilon_0)$ for $Y_0$ around the meeting point of $Z_0$ and $E$. 
\item Milnor balls $B_1$, ..., $B_k$  centered at each of the singular points of $E^{red}$ 

\item tubular neighbourhoods $T_1$, ...$T_r$ contained in $\tilde{X}\setminus \cup_{j=0}^kB_j$ around  each ${E}_i\setminus \cup_{j=0,..,k} B_j$
such that there exist strong deformation retracts $$\zeta_i:T_i\to E_i\setminus \cup_{j=0,..,k} B_j.$$ (see \cite{annals} for technical details)
\end{itemize}
For $s$ small enough we have $$Y_s\subset B_0\cup \bigcup_{j=1}^kB_j\cup \bigcup_{i=0}^rT_i.$$
By the choice of the milnor balls, we have that for $s$ small enough we have trasversality of $Y_s$ and the boundaries of the $B_j$'s and $ T_j$'s 
(see \cite{annals} for technical details).

We are going to give an estimate for $\chi(\calU_{s})$ splitting $\calU_{s}$ as the union of $n^{-1}(B_j)$ and $n^{-1}(T_i)$. Note that $n^{-1}(B_j)$ 
and $n^{-1}(T_j)$ are respectively the normalization of $Y_s\cap B_j$ and $Y_s\cap T_j$. In particular they are dsijoint unions of riemann surfaces with 
boundary. Since we know that $\calU_s$ is a disk and the $B_i$ has transversal boundary with $Y_s$ we have a decomposition of $\calU_{s}$ as in 
Figure \ref{fig:disk}. Furthermore: 

\begin{equation}\label{eq:todo}
\chi(\calU_s)=\chi(n^{-1}(B_0))+\sum_{j=1}^k\chi(n^{-1}(B_j))+\sum_{i=0}^r\chi(n^{-1}(T_j))
\end{equation}
We will separately give estimates for each of the summands in the right hand-side of the equality.

\subsection{Bound in $B_0$}

\begin{figure}\label{fig:B_0}\begin{center}
      \includegraphics[width=40mm]{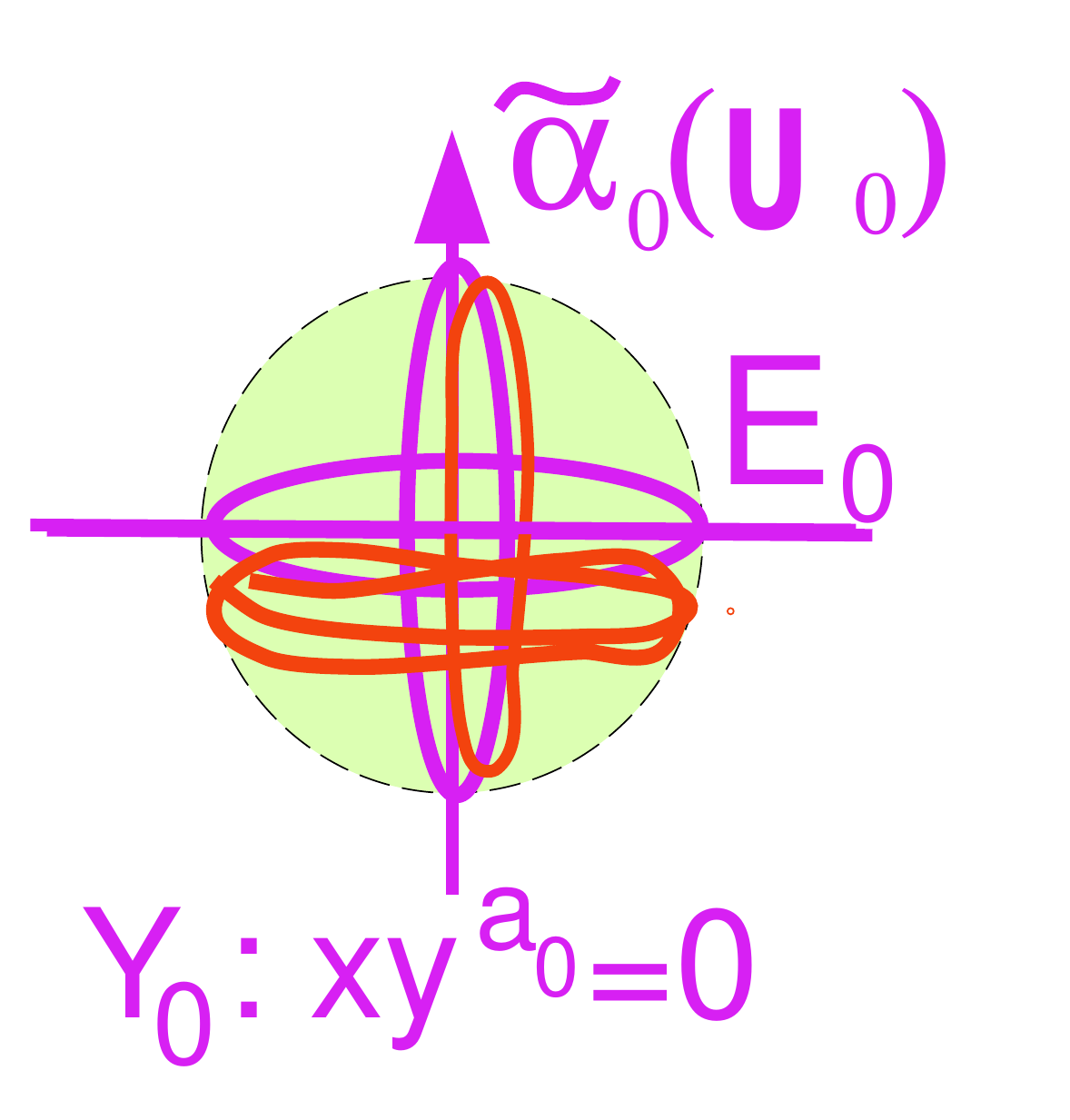}
\caption{Counting the maximal number of disc images in $Y_s\cap B_0$ as a reduced deformation of $Y_0\cap B_0$ of equation $xy^{a_0}=0$ }	
			\end{center}	
\end{figure}

The set $Y_0\cap B_0$ is defined by $f_0(x,y)=xy^{a_0}=0$, where $x$ and $y$ are the coordinates of $B_0$. The divisor $Y_s\cap B_0$ is defined by some deformation $f_s(x,y)=0$ where $f_s$ is a $1$-parameter holomorphic deformation of $f_0$
such that $f_s$ is reduced for $s\neq 0$. Note that 

We observe that $n^{-1}(B_0)$ is a disjoint union of riemann surfaces with boundary. The only connected orientable surface with boundary which has positive Euler characteristic is the disk. Hence $\chi(\calU_s)$ is bounded
above by the number of connected components of $n^{-1}(B_0)$ that are disks. 

There are at the most as many disks in $n^{-1}(B_0)$ as boundary components in $n^{-1}(B_0)$ which are at the most $a_0+1$ since they degenerate to the boundary components of 
$\{xy^{a_0}=0\}\cap B_0$. But by (c) the component of $n^{-1}(B_0)$ whose boundary degenerates to the boundary of $x=0$ in $B_0$, say $A$, is in fact the exterior component of $\calU_s$ (see the Figure \ref{fig:disk}) which can not be a disk, unless it is the whole disk (and this is not possible because $Y_s$ goes outside $B_0$ and meets $E_j$ by (d)). Then $A$ has more than one boundary component. Then, $n^{-1}(B_0)$ has at the most $a_0-1$ disks and we have
\begin{equation}\label{eq:B_0}\chi(n^{-1}(B_0))\leq \# disks \leq a_0-1.\end{equation}

\subsection{Bound in the balls $B_i$, $i\neq 0$}\label{sec:nash_returns}

\begin{figure}\label{fig:B_i}\begin{center}
      \includegraphics[width=40mm]{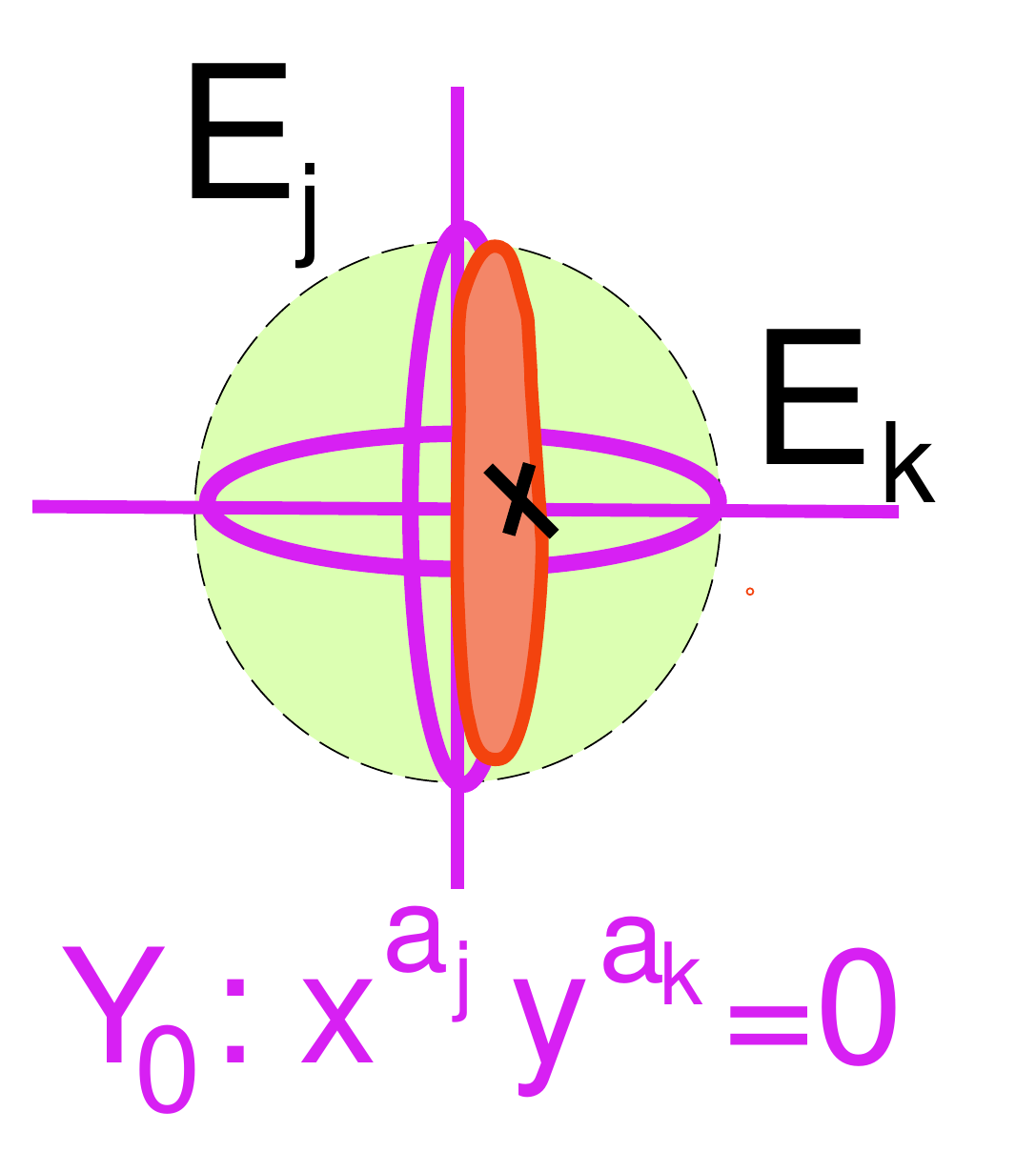}
\caption{Counting the maximal number of disc images in $Y_s\cap B_i$ as a reduced deformation of $Y_0\cap B_i$ of equation $x^{a_j}y^{a_k}=0$ }		
			\end{center}	
\end{figure}

We have 
$Y_0\cap B_i$ defined by $f_0(x,y)=x^{a_j}y^{a_k}=0$, where $x$ and $y$ are the coordinates in $B_i$.

Again, the euler characteristic of 
$n^{-1}(B_j)$ which is a disjoint union of riemann surfaces with boundary is bounded from above by the number of disks. Whenever there is a disk $D$ in $n^{-1}(B_i)$, since its boundary degenerates either to the boundary of $x^{a_j}=0$ or $y^{a_k}=0$ in $B_i$, we have that $n(D)$ will meet at least once either $y=0$ or $x=0$ in $B_0$. Then, 
  
$$\chi(n^{-1}(B_i))\leq \# Y_s\cap Y_0\cap B_i \leq  \sum_{p\in B_i}I_p(Y_s,E).$$

Summing up the estimates of all the balls $B_i$ with $i\neq 0$ we have 
$$\sum_{i=1}^{k}\chi(n^{-1}(B_i))\leq Y_s\cdot E=\sum_k Y_s\cdot E_k.$$
Note that $Y_s\cdot E$ counts the returns (see Section \ref{sec:returns}) with multiplicity. 

Now, we can use that the intersection multiplicity is stable by deformation, that is, $Y_s\cdot E=Y_0\cdot E$ to get
\begin{equation}\label{eq:B_i}\sum_{i=1}^{k}\chi(n^{-1}(B_i))\leq Y_s\cdot E=Y_0\cdot E=(Z_0+\sum_{i=0}^ra_iE_i)\cdot E=1+\sum_{i,k=0,..,r}a_iE_i\cdot E_k.
\end{equation}



\subsection{Bound in every $T_i$}

\begin{figure}\label{fig:T_i}\begin{center}
      \includegraphics[width=40mm]{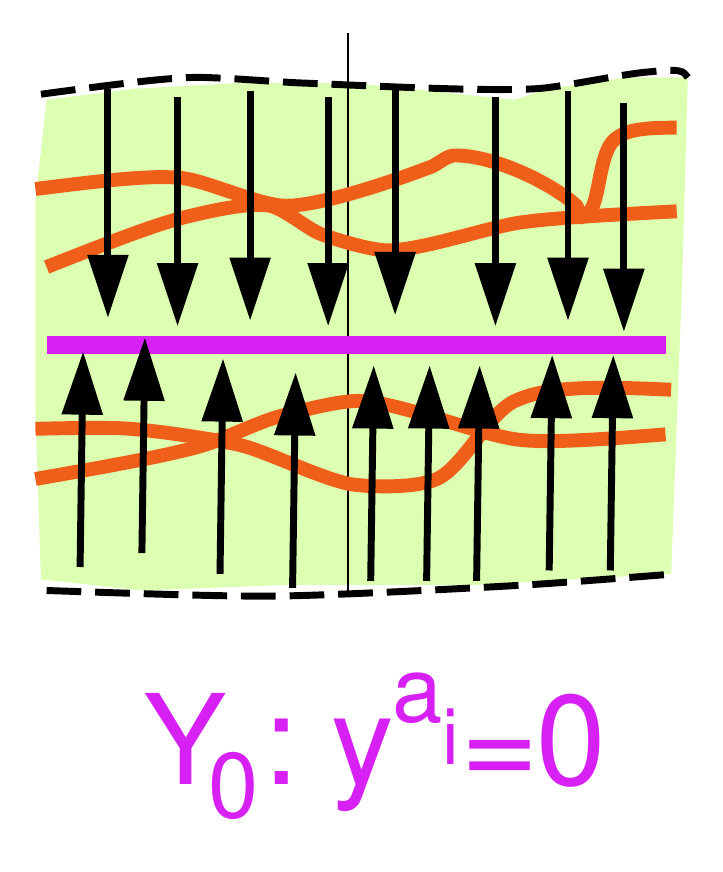}
      \caption{Bounding the euler characteristic of the normalization of $Y_s\cap T_i$ as a reduced deformation of $Y_0\cap T_i$ of equation $y^{a_i}=0$. }
			\end{center}	
\end{figure}

To estimate $\chi(n^{-1}(Y_s\cap T_i))$ we consider the composition
$$\zeta_i\comp n:n^{-1}(Y_s\cap T_i)\to E_i\setminus \cup_jB_j.$$
Although there are some technicalities to be taken into account, for Euler characteristic computations, one may think that it is a holomorphic branched 
cover of Riemann surfaces of degree $a_i$  (the reader may find in~\cite{annals} a completely detailed proof of this).  
Then, by the Riemann-Hurwitz formula we get: 
\begin{equation}
\label{eq:ram1}
\chi(n^{-1}(T_i))\leq a_i\chi(E_i\setminus\cup_j B_j)
\end{equation}

We denote by $g_i$ the genus of $E_i$. For $i\neq 0$, it is clear that $\sum_{k\neq j} E_i\cdot E_k$  counts the number of boundary components of
$E_i\setminus\cup_j B_j$. The number of boundary components  of $E_0\setminus\cup B_i$ is $1+\sum_{k\neq 0} E_0\cdot E_k$ since $E_0$ meets also $B_0$.  

Then, summing up the estimates (\ref{eq:T_i}) for all $i=0,...,r$ we have 
\begin{equation}
\label{eq:T_i}
\sum_{i=0}^r\chi(n^{-1}(T_i))\leq a_0(2-2g_0-1-\sum_{k\neq 0}E_0\cdot E_k)+\sum_{j \neq 0}a_i(2-2g_i-\sum_{k\neq i}E_i\cdot E_k).
\end{equation}

\subsection{Final estimate}

Putting in (\ref{eq:todo}) the estimates (\ref{eq:B_0}), (\ref{eq:B_i}) and (\ref{eq:T_i}) we get that 
\begin{equation}\label{eq:final}
\chi(\calU_s)\leq \sum_i a_i(2-2g_i+E_i\cdot E_i).\end{equation}
 
 By negative definiteness, for any $0\leq i\leq r$, the self-intersection $E_i \cdot E_i$ is a negative integer.
Observe that, since $\pi:\tilde{X}\to X$ is the minimal resolution, for any $0\leq i\leq r$, if $E_i \cdot E_i$ is equal to
$-1$, then either the divisor $E_i$ is singular or it has positive genus 
(otherwise it is a smooth rational divisor with self-intersection equal to
$-1$ and the resolution is non-minimal by the Castelnuovo contractibility Criterion). Since we are assuming that every $E_i$ is smooth we get that 
$2-2g_i+E_i\cdot E_i\leq 0$ for all $i=0,...,r$. 
 
\begin{remark}
\label{re:concluding}

Note that the right hand side is the adjuntion formula in the surface case, which computes the degree of the relative canonical 
 sheaf at each irreducible component of the exceptional divisor. This serves as an inspiration for the higher
dimensional proof of de Fernex and Docampo~\cite{deFDo}.
\end{remark}

Suppose that $\alpha$ is a wedge that does not lift to the minimal resolution. This is equivalent to the existance of indeterminacy of the mapping $\pi^{-1}\circ\alpha$. This implies the inequality $a_0>0$ (which is fact is equivalent) and  bound (\ref{eq:B_0}). 
 
On the other hand, the reader may observe that if the wedge $\alpha$ lifts to the minimal resolution, then the arguments leading to the estimate (\ref{eq:B_0}) break down.
 
The rest of the estimates appearing in our proof are valid in complete generality. So we conclude 
 
\begin{remark}
Our proof shows that if $\alpha$ is a wedge such that the lifting to the minimal resolution of the special arc meets the exceptional divisor in a transverse way, then there is a lifting of $\alpha$ to the minimal resolution.
\end{remark}

\section{The returns of a wedged and deformation theoretic ideas}\label{sec:returns}

Before the proof of the Nash conjecture in \cite{annals}, the second author proved in her PhD the conjecture for the quotient surface singularities (\cite{tesis}, \cite{Pe}). In that proof it is shown that, despite the local nature of arcs, at least semi-local techniques were needed in order to study the arc space and the existence of certain families of arcs or wedges. 

In particular, it was observed that for a representative $\alpha|_U$ of a wedge, returns might be non-avoidable. We recall that a \emph{return} is a point in $\alpha_s^{-1}(0)$ different from the origin for $s\neq 0$. A return $p\in \alpha_s^{-1}(0)$ is identified with the associated arc that consists in viewing $\alpha_s|_{U_s}$ as a germ at $p$.   
 If one thinks of an arc as the image of a parametrization ``starting at the singular point'', the returns are the points where the 
parametrization passes again through the singular point. The contribution of the returns is crucial in the Euler characteristic estimates needed in our proof of the Nash conjecture (see Section \ref{sec:nash_returns}).

The study of returns had a direct impact in the application of the valuative criterion to rule out adjacencies $N_E\subset N_F$ in \cite{Pe}. The valuative criterion was first studied in \cite{Re1}, \cite{Pl}. 
Given an exceptional prime divisor $D$ over $(X,0)$, we denote by $ord_D$ the associated divisorial valuation. The valuative criterion says that if there exists a germ $g\in \calO_X$ such that $ord_E(g)<ord_F(g)$ where $E$ and $F$ are exceptional prime divisors, then the adjacency $N_E\subset N_F$ is not possible. Now, taking into account the retunrs, we can say that  if we have an inequality $ord_E(g)<ord_F(g)+ord_{F'}(g)$ then there is no wedge realising the adjacency  $N_E\subset N_F$ with a return with lifting  by $F'$ (nor a wedge realising the adjacency  $N_E\subset N_{F'}$ with a return by $F$). 

This idea is applied  in \cite{Pe} more conveniently for the pullback of the wedges by the quotient map $q:(\CC^2,0)\to (X,0)$ for a quotient surface singularity. 

In \cite{annals} this criterion is completely understood in Section 3.2. as follows.   
Consider a wedge realising an adjacency $N_{E_0}\subset N_{E_j}$ as in Section \ref{representatives}.
Since the divisor $Y_s$ defined in (\ref{pijotada})-(\ref{eq:Ys}) is a deformation of the divisor $Y_0$, we have the equality 
\begin{equation}
\label{eqns}Y_0\centerdot E_i=Y_s\centerdot E_i\end{equation}
for any $i$. Recall notation $Y_0=Z_0+\sum_ia_iE_i$  . Denote by $b_i$ the intersection product of $Y_s\centerdot E_i$ and by $M$ the matrix of the intersection
form in $H_2(\tilde{X}_\epsilon,\ZZ)$ with respect to the basis $\{[E_0],...,[E_r]\}$.
Then, (\ref{eqns}) can be expressed as follows:
\begin{equation}
\label{eq:sistemalineal}
M(a_0,...,a_r)^t=(1-b_0,b_1,...,b_r)^t.
\end{equation}
The number $b_i$ is the number
of \emph{returns} of the wedge through the divisor $E_i$ counted with appropiate multiplicity.


An important observation is that all the entries of the inverse matrix $M^{-1}$ are non-positive (see Lemma 10 in \cite{annals}).

The equality (\ref{eq:sistemalineal}) can be used to prove that wedges realizing certain adjacencies with certain prescribed returns $b_i$ do not exist:   
 the existence of such a wedge is impossible if the solution $a_0$,...,$a_n$ of (\ref{eq:sistemalineal}) has either a negative or a non-integral entry.
 
Moreover, to finish the proof in \cite{Pe} for the $E_8$ singularities, further arguments using deformation theory were needed. 
There, wedges realising an adjacency $N_{E_0}\subset N_{E_j}$ with a given special arc are seen as $\delta$-constant deformations of the cuve 
parametrized by the special arc. Then, the versal deformation of the  curve parametrized by the special arc is computed. The codimensions of the 
$\delta$-constant stratum and the codimension of the stratum of curves with the  topological type of the generic curve of a family parametrized by a 
wedge representative with  prescribed returns are computed. 
The inequality that these codimensions satisfy is not compatible with the existance of such a wedge (see Proposition 4.5. in \cite{Pe}).

\section{The proof by de Fernex and Docampo for the higher dimensional case}

De Fernex and Docampo figured out an 
algebro-geometric proof of the Nash conjecture based on bounds of coefficients of suitable relative canonical sheaves~\cite{deFDo}.
This enabled them to formulate and prove a correct statement of the Nash correspondence in higher 
dimension. 

A \emph{terminal valuation} is a 
divisorial valuation on $X$ such that there exists a terminal minimal model $\pi:Y\to X$ of $X$ such that the center of the valuation 
is a divisor in $Y$. The centers of terminal valuations are essential divisors. The main theorem of de Fernex and Docampo is

\begin{theo}[de Fernex, Docampo]
Terminal valuations are at the image of the Nash map.
\end{theo}
 
The begining of the proof is similar to the surface case: assuming that the result is false they derive the existence of a wedge 
such that its special arc lifts to $Y$ in a transverse way to the center of a terminal valuation, but that can not be lifted to $Y$.
Afterwards, assuming the existence of such a wedge they bound a coefficient for a relative canonical sheaf in two different ways and produce a contradiction.
So, in fact they prove:

\begin{theo}[de Fernex, Docampo]
\label{th:dd}
Let $\pi:Y\to X$ be a terminal model of $X$. Any wedge $\alpha$ such that its special arc lifts to $Y$ in a transverse way to the center 
of a terminal valuation, admits a lifting to $Y$.
\end{theo}

We refer to the original paper for a complete explanation of their proof. 
Here, instead we explain the main ideas of the proof in the context of surfaces. In this case $Y=\tilde{X}$ where $\tilde{X}$  is the unique terminal model which is
the minimal resolution $\pi:\tilde{X}\to X$.

On one hand it is easier to digest, and all main ideas
appear in this case. On the other hand, by doing it we derive a precise set of numerical equalities (see (\ref{eq:numerical})) that are satisfied for 
any non-constant wedge $\alpha:(\CC^2,O)\to X$ not lifting to the minimal resolution and such that $\pi^{-1}\comp\alpha$ is a meromorphic map with an only indeterminacy point, 
regardless this wedge has special arc lifting transversely or not (this means any wedge that is used in practice). 
This set of equalities have not been observed before. If one assumes that the special arc of the wedge lifts transversely, one may derive
a chain of inequalities giving the contradiction in a straightforward way from this set of equalities.

Let $\alpha:(\CC^2,O)\to X$ be any wedge so that $\pi^{-1}\comp\alpha$ is a rational map from $(\CC^2,O)$ to $\tilde{X}$,
and such that the special arc of the wedge .
Let $\sigma:Z\to (\CC^2,O)$ be the minimal sequence of blow ups at points resolving the 
indetermination of $\pi^{-1}\comp\alpha$. Let $\beta:Z\to\tilde{X}$ be the map such that $\pi\comp\beta=\alpha\comp\sigma$. De Fernex and
Docampo shift the computation from $\tilde{X}$ to $Z$. 

Let $F=\sum_{i=1}^mF_i$ and $E=\sum_{i=1}^nE_i$ be the decomposition in irreducible components of $\sigma$ and $\pi$ respectively. 
In $(\CC^2,O)$ we consider coordinates $(t,s)$ so that $t$ is the arc variable and $s$ is the deformation parameter. The special 
arc of the wedge is then $\alpha(t,0)$. We order the components so that $F_1$ is the unique component where the strict transform
of $V(s)$ meets.

Denote by $K_Z=\sum_{i}a_iF_i$ the canonical divisor of $Z$. It is the only representative of the canonical class $K_Z$ supported 
at the exceptional divisor. Each $a_i$ is positive, and a simple observation on the behaviour of the 
canonical divisor under blow up shows that 
\begin{remark}
\label{re:a1}
The number $a_1$ is the number of blowing up centers touching the strict transform of $V(s)$.
\end{remark}

Since $\beta$ is a morphism between smooth spaces, the relative canonical class $K_{Z/\tilde{X}}:=K_Z-\beta^*K_{\tilde{X}}$ may be
represented by the divisor associated with the jacobian of $\beta$. This is an effective divisor. When we write $K_{Z/\tilde{X}}$ 
we mean such a divisor. We decompose it as 
$$K_{Z/\tilde{X}}=K_{Z/\tilde{X}}^{exc}+K_{Z/\tilde{X}}^{hor},$$
where $K_{Z/\tilde{X}}^{exc}$ is the part with support on the exceptional divisor of $\pi$ and $K_{Z/\tilde{X}}^{hor}$ the complement.
We have the equality
\begin{equation}
\label{eq:clave}
K_Z-K_{Z/\tilde{X}}^{exc}=K_{Z/\tilde{X}}^{hor}+\beta^*K_{\tilde{X}}.
\end{equation}
The left hand side is a divisor concentrated in the exceptional set of $\sigma$. 

In order to express the right hand side as a divisor
concentrated in the exceptional set we let $M$ be the intersection matrix of the collection of divisors $\{F_i\}$ in $Z$. Since $\sigma$
is a sequence of blow ups we have that $M$ is unimodular and that its inverse $M^{-1}$ is the matrix whose $i$-th column 
$(m_{1,i},...,m_{n,i})^t$ is obtained as follows: let the $C_i$ be the curve in $(\CC^2,O)$ given by the image of $\sigma$ of a curvette
transverse to $F_i$. Then its total transform to $Z$ is
$$\sigma^*C_i=\sum_j-m_{j,i}F_j.$$
As a consequence we obtain that all the entries of $M^{-1}$ are strictly negative (this is a general phenomenon which is well known,
see for example Lemma 10 in \cite{annals} for the proof for general normal surface singularities). 

If we express $K_{Z/\tilde{X}}^{exc}=\sum_i b_iF_i$, and define the 
intersection numbers
$$c_i:=K_{Z/\tilde{X}}^{hor}\centerdot F_i,$$
$$d_i:=\beta^*K_{\tilde{X}}\centerdot F_i,$$
Equality (\ref{eq:clave}) becomes the following system of numerical 
equalities
\begin{equation}
\label{eq:numerical}
(a_1-b_1,...,a_n-b_n)^t=M^{-1}(c_1+d_1,...,c_n+d_n).
\end{equation}

This numerical equality works for any resolution $\pi:\tilde{X}\to X$ (non-necessarily minimal) and any wedge $\alpha$, so it may 
have applications in other problems. One instance could be the generalized Nash problem explained in Section \ref{sec:NashGen}.

Observe that, since $K_{Z/\tilde{X}}^{hor}$ has no component included in the exceptional divisor, each $c_i$ is non-negative.

If we asume now that $\pi$ is the minimal resolution, we have
$$d_i=\beta^*K_{\tilde{X}}\centerdot F_i=K_{\tilde{X}}\centerdot \beta_*(F_i),$$
which is non-negative by adjunction formula, using the fact that $\pi:\tilde{X}\to X$ is the minimal resolution.

This means that the right hand side in Equation~(\ref{eq:numerical}) is {\em non-positive}.

In order to prove Theorem~\ref{th:dd} for the surface case we assume that the wedge has the special arc lifting transversely to the
exceptional divisor, and estimate the coefficient $a_1-b_1$ in the left hand side of Equation~(\ref{eq:numerical}).

By Remark~\ref{re:a1} if the wedge $\alpha$ does not lift to $\tilde{X}$ then $a_1$ is strictly positive and integral. 
Since the right hand side of 
Equation~(\ref{eq:numerical}) is non-positive, in order to finish the proof it is enough to prove the strict inequality
$$b_1<1.$$

Let $\rho:Z\to Z'$ and $\beta': Z'\to \tilde{X}$ be such that  the factorization $\beta=\beta'\comp\rho$ consists in collapsing all non-dicritical 
components of the exceptional divisor in $Z$ (a component $F_i$ is non-dicritical if $\beta(F_i)$ is a point).
The surface $Z'$ has sandwiched singularities, which are rational and $\QQ$-Gorenstein. Then the canonical divisor 
$K_{Z'}$ is $\QQ$-Cartier. Therefore we may define the relative canonical class $K_{Z/Z'}$. This class has a 
unique representative as a $\QQ$-divisor supported in the exceptional set of $\rho$ such that all its coefficients are non-positive.

We have the equality $K_{Z/\tilde{X}}=K_{Z/Z'}+\rho^*K_{Z'/\tilde{X}}$. By the non-positivity of the coefficients of $K_{Z/Z'}$ we get
$$b_1=ord_{F_1}(K_{Z/\tilde{X}})=ord_{F_1}(K_{Z/Z'}+\rho^*K_{Z'/\tilde{X}})\leq ord_{F_1}(\rho^*K_{Z'/\tilde{X}}).$$

We make an abuse of language and denote the components of the exceptional divisor of $Z'$ by the same name that they have in $Z$.
In order to estimate $b_1$ we enumerate $\{F_{i_1},...,F_{i_l}\}$ the components of the exceptional set of $Z'$ which 
contain the image by $\rho$ of $F_1$. This is a subset of the components of $F$ not collapsed by $\rho$. Observe that if $F_1$ is 
dicritical then this set of components has $F_1$ as a unique element, and the estimate that we will prove right away becomes much easier.  

The special arc of the wedge $\alpha$ lifts transversally through an irreducible component of $E$. We enumerate the components so that 
this component is $E_1$. Then for each of the components $F_{i_j}$ we have that $\beta(F_{i_j})=E_1$. 

Before proving our final estimate we need the following observation: the following equality holds
$$ord_{F_{i_j}}(K_{Z'/\tilde{X}})=ord_{F_{i_j}}((\beta')^*E_1)-1.$$
This holds because $K_{Z'/\tilde{X}}$ is given at smooth points by the divisor associated with the jacobian of $\beta'$, and at a generic point 
of $F_{i_j}$ the mapping $\beta$ can be expressed in local coordinates as $\beta(u,v)=(u^a,v)$, where $a=ord_{F_{i_j}}((\beta')^*E_1)$.

The last estimate we need is:
$$ord_{F_1}(\rho^*K_{Z'/Y})= ord_{F_1}(\sum_{j=1}^l ord_{F_{i_j}}(K_{Z'/\tilde{X}})\rho^*F_{i_j}))=$$
$$=\sum_{j=1}^l ord_{F_1}(ord_{F_{i_j}}((\beta')^*E_1)-1)\rho^*F_{i_j})<\sum_{j=1}^l ord_{F_1}(ord_{F_{i_j}}((\beta')^*E_1))\rho^*F_{i_j})=$$
$$=ord_{F_1}\rho^*(\beta')^*E_1=ord_{F_1}\beta^*E_1=1.$$

The last equality holds because the special arc lifts transversely by $E_1$. This concludes the proof.

\section{On the generalized Nash problem and the classical adjacency problem}\label{sec:NashGen}

Let $X$ be a normal surface singularity. \emph{The Generalized Nash problem} consists in characterizing the pair of divisors $E$, $F$ appearing in resolutions of $X$ such that 
the adjacency $N_F\subset N_E$ holds. 

For our proof of the Nash conjecture it is essential to construct a holomorphic wedge $\alpha$, as we have explained before. This was achieved in~\cite{Bo}. The technique developed to achive this
gave, as a by product, a proof of the fact that the validity of the Nash conjecture only depends on the topology of the link of the surface singularity, or equivalently, in the combinatorics
of the minimal good resolution. The same technique could be adapted to prove that the generalized Nash problem is a topological problem in the following sense.

Since the generalized Nash problem is wide open even in the case in which $X$ is smooth we concentrate in this case. To any $2$ exceptional divisors $E$ and $F$ of a sequence of blow ups at the 
origin of $X$, we may associate a decorated graph as follows: consider the minimal sequence of blow ups of $\pi:Y\to X$ where both  $E$ and $F$ appear. Decorate the dual graph of the 
exceptional divisor of $\pi$ attaching to each vertex the weight given by the self-intersection of the corresponding divisor. Finally add labels $E$ and $F$ to the vertices corresponding to
the divisors $E$ and $F$ respectively. In~\cite{BoPP} we proved:

\begin{theo}
Let $(E_1,F_1)$ and $(E_2,F_2)$ be two pairs of divisors having the same associated graph. Then the adyacency $N_{F_1}\subset N_{E_1}$ is satisfied if and only if the adyacency 
$N_{F_2}\subset N_{E_2}$ is satisfied.
\end{theo}

As a consequence we could improve the discrepancy obstruction for adyacencies: see Corollary 4.17 and 4.19 of \cite{BoPP}. Furthermore, we get a nice structure of nested nash sets in the arc space of $\CC^2$ and we made some conjectures about it: see conjectures 1 and 2 of \cite{BoPP}.

For the sake of completeness we summarize very briefly the other main result of \cite{BoPP}. 

Given a prime divisor $E$ over the origin of $X$ we consider its associated valuation $\nu_E$. It is easy to see that if we have the adjacency  $N_F\subset N_E$ then the inequality
$\nu_E\leq \nu_F$ holds. However this criterion is not enough to characterize the Nash adjacencies (see Section \ref{sec:returns}). Our second main result is a characterization of the previous inequality in terms of deformations 
of plane curves. 

We say that a plane curve germ $C$ is \emph{associated with} $E$ in a model $\pi:Y\to \CC^2$ where $E$ appears if its strict transform by $\pi$ meets $E$ trasversely at a point which does not meet 
the singular set of the exceptional divisor of $\pi$. A deformation $g_s$ of function germs is a holomorphic function depending holomorphically on a parameter $s$. It is linear if it is of 
the form $g_0+sh$ for $g_0$ and $h$ holomorphic.    

\begin{theo}
Let $E$ and $F$ be prime divisors over the origin of $\CC^2$ and let $S$ be the minimal
model containing both divisors. The following are equivalent:
\begin{enumerate}
 \item $\nu_E\leq \nu_F$,
\item There exists a deformation $g_s$ with $g_0$ associated with $F$ in $S$ and $g_s$
 associated with $E$ in $S$, for $s\neq 0$ small enough,
\item There exists a linear deformation $g_s$ with $g_0$ associated with $F$ in $S$ and $g_s$
 associated with $E$ in $S$, for $s\neq 0$ small enough.
\end{enumerate}
\end{theo}

In fact in~\cite{BoPP} we prove a more general version which allows non-prime divisors $E$ and $F$.

This Theorem provides a very easy way of producing adjacencies of plane curve singularities. Using it we were able to recover most Arnol'd adjacencies. See Section~3.4 of \cite{BoPP} for 
detailed explanations.

\section{Holomorphic arcs}

In the proof of the Nash conjecture for surfaces we study arcs and wedges from a convergent viewpoint, and take representatives. In this sense a wedge is for us a deformation of holomorphic
maps from a disc to a representative of the singularity. At a generic parameter the preimage of the singular point in general contains more points in the disc that just the origin. These
points are unavoidable and we call them {\em returns} following \cite{Pe} (see Section \ref{sec:returns}). 

J. Koll\'ar and A. Nemethi started in~\cite{KoNe} the sytematic study of convergent arc spaces as oposed to the classical formal arc spaces. We briefly summarize their main results and 
questions.

Let $D$ denote the closed unit disc. A holomorphic map defined on $D$ is the restriction to $D$ of a holomorphic map in an open neighbourhood of $D$. Let $X$ be a singularity.

\begin{definition}
A {\em complex analytic} arc is a holomorphic map $\gamma:D\to X$ such that the preimage of the singular set does not intersect the boundary $\partial D$. A {\em short complex analytic arc}
is a complex analytic arc such that the preimage of the singular set is just $1$ point. A {\em deformation of a complex analytic arc} parametrized by an analytic space $\Lambda$ 
is a holomorphic map 
$$\alpha:D\times\Lambda\to X$$
such that each for any $s\in\Lambda$ the restriction $\alpha_s:=\alpha|_{D\times\{s\}}$ is a complex analytic arc. 
If all the arcs appearing are short complex analytic arcs we say that $\alpha$ is a {\em deformation of short complex analytic arcs}.
\end{definition}

\begin{remark}
What we do in Section~\ref{representatives} is to derive from a convergent wedge a deformation of a complex analytic arc parametrized by a disc $\Lambda$. The special arc at this deformation
is a short arc, but the generic arcs $\alpha_s$ are not short arcs in general, due to the existence of returns. The topological analysis of that deformation of complex analytic arcs 
yields the proof of the Nash conjecture.
\end{remark}

Denote by $Arc(X)$ and by $ShArc(X)$ the sets of convergent analytic arcs in $X$. Koll\'ar and Nemethi give natural metrics on these spaces, which endow them with a topology. 
Given an arc $\gamma:D\to X$, since the preimage of the singular set $Sing(X)$ is disjoint from the circle $\partial D$, the restriction $\gamma|_{\partial D}$ defines an element of 
the fundamental group modulo conjugation $\pi_1(X\setminus\Sing(X))/(conjugation)$. 
Since this element does not change by continuous deformation of the arc $\gamma$ we have defined ``winding number maps''
$$\pi_0(Arc(X))\to \pi_1(X\setminus\Sing(X))/(conjugation),$$
$$\pi_0(ShArc(X))\to \pi_1(X\setminus\Sing(X))/(conjugation).$$

The main result in~\cite{KoNe} concerns short arcs:

\begin{theo}[Koll\'ar, Nemethi]
The winding number map 
$$\pi_0(ShArc(X))\to \pi_1(X\setminus\Sing(X))/(conjugation)$$
is injective for any normal surface singularity $X$. It is bijective for quotient surface singularities. 
\end{theo}

For general normal surface singularities the winding number map is far from being surjective, but its image is described in~\cite{KoNe} in terms of the combinatorics of the 
resolution, or what is the same, the topology of the link.

On the other hand, the winding number map for long arcs, which is the one that is 
more related with the original Nash question, is not well understood.
\begin{problem}[Koll\'ar, Nemethi]
Is the winding number map
$$\pi_0(Arc(X))\to \pi_1(X\setminus\Sing(X))/(conjugation)$$
injective?
\end{problem}

In~\cite{KoNe} many other open problems are proposed. Some interesting ones are concerned with the definition of a ``finite type'' holomorphic atlas in $Arc(X)$ 
(see Conjecture 72 of loc. cit.), and with the existence of a curve selection Lemma in $Arc(X)$.

\end{document}